\documentclass[12pt]{amsart}
\usepackage {amsmath, amscd, amsbsy, amsfonts}
\usepackage{amsmath,latexsym,amssymb,amsfonts,mathrsfs,bbm}
\usepackage{fullpage}
\usepackage{graphicx}
\usepackage{hyperref}
\newlength{\baseunit}               
\newcommand{\R}{\mathbbmss{R}}

\newcommand{\s}{\mathbf{s}}

\def \Y {{\mathcal Y}}
\def \gam {{\Gamma}}

\def \n {{\mathcal N}}

\def \D {\mathcal D}
\def \Q {\mathbb Q}
\def \proj {{\mathbb{P}^r}}

\def \X {\mathcal X}

\def \d {\Delta}
\def \codim {\text {codim}}

\def \L {{\mathcal L}}
\def \W {{\mathcal W}}

\def \K {{\mathcal K}}

\def \F {{\mathcal F}}

\def \s {{\mathcal S}}
\def \R {{\mathcal R}}
\def \C {{\mathcal C}}

\def \n {{\mathcal N}}
\def \mbar {\overline{\mathcal M}}
\def \sep {\ || \ }
\def \wt {\widetilde}
\def \E {{\mathcal E}}

\newtheorem{theorem}{Theorem}[section] 

\newtheorem{lemma}[theorem]{Lemma}

\newtheorem{corollary}[theorem]{Corollary}

\newtheorem{proposition}[theorem]{Proposition} 
 
\newcommand{\bpf}{ {\it Proof.  }}
\newcommand{\epf}{\qed \vspace{+10pt}}

\title[Enumerative Invariants.]{Characteristic numbers of rational cuspidal space curves}

\author{Dung Nguyen}

\begin{document}
\maketitle
\begin{abstract}
We solve the problem of computing characteristic numbers
of rational space curves with a cusp, where there may or 
may not be a condition on the node. The solution is given 
in the form of effective recursions. We give explicit formulas when the dimension
of the ambient projective space is at most $5$. Many
numerical examples are provided. A C++ program
implementing most of the recursions 
is available upon request.
 \vspace{-20pt}
\end{abstract} 
\section{introduction}
 Charateristic number of curves in projective spaces is a classical problem in algebraic
geometry: how many curves in $\proj$ of given degree and genus 
that pass through a general set
of linear subspaces, and are tangent to a general set of hyperplanes? 
The advent of stable maps has provided a powerful tool to attack
such problem, especially in low genus. In this paper,
we use stable maps to solve the problem of counting
rational cuspidal space curves, where one may impose
tangency conditions and an incident condition on the
cusp. This is a quick application of the results in
\cite{dn}.

 The number of rational cuspidal plane curves satisfying
incidence conditions were compute in \cite{p}. The method
was to use intersections of divisors on $\mbar_{0,n}(2,d)$,
as the locus of cuspidal plane curves is a divisor on
$\mbar_{0,n}(2,d)$. This is no longer true for cuspidal
space curves. The incidence-only characteristic numbers
of rational cuspidal space curves were computed in \cite{zg}.
There are a number of classical results regarding full
characteristic numbers of some family of rational
cuspdial curves. The numbers of cuspidal plane cubics
satisfying incidence and tangency conditions, where
the cusp may or may not be subject to an incident
condition, were computed in \cite{luf}. The analogous
numbers for cuspidal plane cubics in $\mathbb{P}^3$ 
were computed in \cite{hm}

 Our approach is simple. We use the well-known
divisorial relation on $\mbar_{0,4}(r,d)$, that
is the pull-back of the trivial relation
on $\overline{M}_{0,4}.$ We intersect this relation
with a carefully choosen substack whose general
members are maps that have a node in the image
(rational nodal curves). Since cuspidal curves
can be viewed as limits of rational nodal curves,
the locus of cuspidal curves form a boundary
component of the substack of rational nodal curves.
Thus the above relation allows us to relate
enumerative invariants of cuspidal curves to
those of rational nodal curves, which can be
computed from \cite{dn}.

 The main advantage of this method is that it
is easy to implement and works well for all type
of conditions. For example, there is no need
to separate the family of rational cuspidal space
curves into subfamilies where the cusp lie on
general linear subspaces of different codimensions, as
all characteristic numbers of such families can
be computed using one common recursion. Both the counting
of rational cuspidal space curves in this note,
and the counting of rational nodal space curves in \cite{dn}
will have application in couting elliptic space curves,
in an upcoming paper by the author.

 The content of the paper is as follows. In Section 2, we introduce the basic
notations that will be used throughout
this note. In Section 3, we propose
and prove the recursion relating
the enumerative geometry of cuspidal
curves to that of rational nodal curves
and that of reducible curves whose
components are smooth rational intersecting
at two points. Section 4 gives numerical
examples and compare them to known results
in literature.

\section{Definitions and Notations}
\subsection{The moduli space of stable maps of genus $0$ in $\proj$.}
As usual, $\mbar_{0,n}(r,d)$ will denote the Kontsevich compactification of the moduli
space of genus zero curves with $n$ marked points of degree $d$ in $\proj$. 
 We will also be using the notation $\mbar_{0,S}(r,d)$ where the markings 
are indexed by a set $S$.
The following are Weil divisors on $\mbar_{0,S}(r,d)$: 
\begin{itemize}
 \item The divisor $(U \sep V)$ of $\mbar_{0,S}(r,d)$
 is the closure in $\mbar_{0,S}(r,d)$ of the locus 
of curves with two components such that $U \cup V = S$ is a partition 
of the marked points over the two components.
\item The divisor $(d_1,d_2)$ is the closure in $\mbar_{0,S}(r,d)$
 of the locus of curves with two components, sucht that $d_1 + d_2 = d$ 
is the degree partition over the two components. 
\item The divisor $(U,d_1 \sep V,d_2)$ 
is the closure in $\mbar_{0,S}(r,d)$ of the locus of 
curves with two components, where $U \cup V = S$
and $d_1 +d_2 = d$ are the partitions of markings 
and degree over the two components respectively.
\end{itemize}

\subsection{The constraints and the ordering of constraints.}

 We will be concerned with the number of curves passing through a constraint. Each 
constraint is denoted by a $(r+1)-$tuple $\d$ as follows : \\
{\bf (i)}  $\d(0)$ is the number of hyperplanes that the curves need to be tangent to. \\
{\bf (ii)} For $0<i \leq r$, $\d(i)$ is the number of subspaces of codimension
$i$ that the curves need to pass through. \\
{\bf (iii)} If the curves in consideration have a node (or cusp) and we place a condition
on the node (cusp), that is the node (cusp) has to belong to a general codimension
$k$ linear subspace, then $\d$ has $r+2$ elements and the last element, $\d(r+1)$, is $k$.

 Note that because in general a curve of degree
$d$ will always intersect a hyperplane at $d$ points, introducing an incident condition
with a hyperplane essentially means multiplying the cycle class cut out
by other conditions by $d$.
For example, if we ask how many genus zero curves of degree $4$ in $\mathbb P^3$ that pass through the constraint
$\d = (1,2,3,4,0)$, that means we ask how many genus zero curves of degree $4$ pass through
three lines, four points, are tangent to one hyperplane, and then multiply that answer
by $4^2$. We will also refer to $\d$ as a set of linear spaces,
hence we can say, pick a space $a$ in $\d$. \\ \\
 We consider the following ordering on the set of constraints, in order to prove that our algorithm will 
terminate later on. Let $r(\d) = -\sum_{i >1}^{i\leq r} \d[i]\cdot i^2$, and this will be our rank function.
 We compare two constraints $\d,\d'$ using the 
following criteria, whose priority are in the following order :
\begin{itemize}
 \item  If $\d(0) = \d'(0)$ and $\d$ has fewer non-hyperplane elements than $\d'$ does, then $\d<\d'$. 
 \item If $\d(0) > \d'(0)$  then $\d < \d'$.
\item If $r(\d) < r(\d')$ then $\d < \d'$. 
\end{itemize}
Informally speaking, characteristic numbers
where the constraints are more spread out at two ends 
are computed first in the recursion.
We write $\d = \d_1\d_2$ if $\d = \d_1 + \d_2$
as a parition of the set of linear spaces in $\d.$

\subsection{The stacks $\R,\n, \R\R, \R\R_2, \s(r,d)$ .}
 We list the following definitions of stacks of stable maps that will
occur in our recursions. \\ \\
{\bf 1)} Let $\R(r,d)$  be the usual moduli space of genus zero stable maps
$\mbar_{0,0}(r,d)$. \\ \\
{\bf 2)} Let $\n(r,d)$ be the closure in $\mbar_{0,\{A,B\} }(r,d)$ of the locus of maps of smooth
rational curves $\gamma$ such that $\gamma(A) = \gamma(B)$. Informally,
$\n(r,d)$ parametrizes degree $d$ rational nodal curves in $\proj$. \\ \\
{\bf 3)} For $d_1,d_2>0,$ let $\R\R(r,d_1,d_2)$ be $\mbar_{0,\{C\}}(r,d_1) \times \mbar_{0,\{C\}}(r,d_2)$
where the fibre product is taken over evaluation maps $ev_{C}$ to $\proj.$ \\ \\
{\bf 4)} Similarly we can define $\n\R(r,d_1,d_2)$, $\E\R(r,d_1,d_2)$ (see figure 1).\\ \\
{\bf 5)} For $d_1,d_2>0$, let $\R\R_2(r,d_1,d_2)$ be the closure in $\mbar_{0,\{A,C\}}(r,d_1) \times_{\proj} \mbar_{0,\{B,C\}}(r,d_2)$  
(the projections are evaluation maps $e_C$) of the locus of maps $\gamma$ 
such that $\gamma(A) = \gamma(B)$. We call maps in this family
rational two-nodal reducible curves.
{\bf 6)} For $d>0$ let $\s(r,d)$ be the closure in $\mbar_{0,\{A\}}(r,d)$ of the locus of 
maps of smooth domains $\gamma$ such that the differential vanishes at $A$. Informally,
$\S(r,d)$ parametrizes degree $d$ rational cuspidal curves in $\proj.$

$$\includegraphics[width = 150mm]{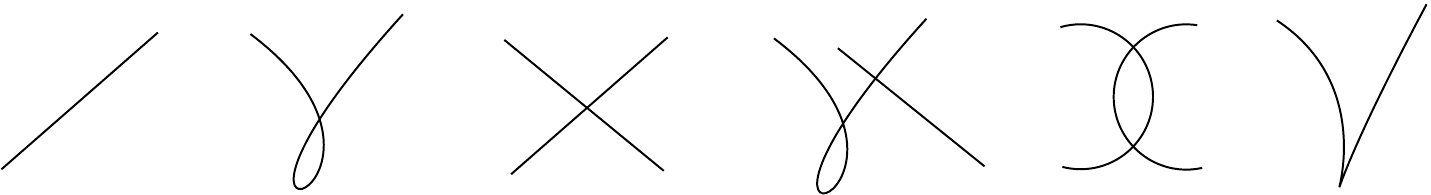}$$ \\
$$\text{Fig 1. Pictorial description of a general curve in the stacks $\R,\n, \R\R, \n\R, \R\R_2, \s(r,d)$}$$

\subsection{Stacks of stable maps with constraints.} Let $\F$
be a stack of stable maps of curves into $\proj$. For a constraint
$\d$, we define
$(\F, \d)$ be the closure in $\F $ of the locus of maps
that satisfy the constraint $\d$. If the stack of maps
$\F$ has two marked points $A$ and $B$, we define
$(\F, \L_A^u\L_B^v)$ to be the closure
in $\F$ of the locus of maps $\gamma$ such 
that $\gamma(A)$ lies on $u$ general hyperplanes,
 and that  $\gamma(B)$ lies on $v$ general hyperplanes. 

If $\F$ has one marked point $A$ then 
we define $(\F,\L_A^u\W^v_A)$ 
to be the closure of maps $\gamma$ such that
$\gamma(A)$ lies on $u$ general hyperplanes, and that
the image of $\gamma$ is smooth at $\gamma(A)$ and the tangent line to the image of $\gamma$
at $\gamma(A)$ passes through $v$ general codimension $2$
subspaces. 

If a stack of $\F$ consists of a finite
number of points then we denote
$\# \F$ to be the stack-theoretic length of $\F$.

If $\F$ is a closed substack of the stacks $\n\R, \R\R$ then we denote
$(\F,\gam_1,\gam_2,k)$  to be the closure in $\F$
of the locus of maps $\gamma$ such that the
restriction of $\gamma$ on the $i-$th component
satisfies constraint $\gam_i$ and that
$\gamma(C)$ lies on $k$ general hyperplanes. We use the notation $(\F,\d,k)$ if we don't want
to distinguish the conditions on each component.

If $\F$ is a closed substack of $\R\R_2(r,d_1,d_2)$ then we denote
$(\F,\gam_1,\gam_2,k,l)$  to be the closure in $\F$
of the locus of maps $\gamma$ such that the
restriction of $\gamma$ on the $i-$th component
satisfies constraint $\gam_i$ and that
$\gamma(C)$ lies on $l$ general hyperplanes, and that
$\gamma(A) = \gamma(B)$ lies on $k$ general hyperplanes.
Similary,
we use the notation $(\F,\d,k,l)$ if we don't want
to distinguish the conditions on each component. 

Note that for maps of reducible source curves, tangency
condition include the case where the image of the node
lies on the tangency hyperplane, as the intersection
multiplicity is $2$ in this case.

\section{Counting Rational Cuspidal Curves via Rational Nodal Curves}
We need a result to establish the locus of
rational cuspidal curves as a boundary component
of the locus of rational nodal curves (both
are substacks of the genus $0$ stable map 
space).
\begin{lemma}
  Let $\n(r,d)$ be as in Section 2. Let $\K$ be the closure
in $\n(r,d)$ of maps $\gamma$ such that : \\
(i) The domain has two components and both $A,B$ are on the same component \\
(ii) $\gamma$ contracts that component. \\
Let $C$ be the node of the source curve. Then if 
$\gamma$ is a general map in $\K$ 
then the restriction of $\gamma$ on the
other component has differential
vanished at $C$.
\end{lemma}
\bpf We can look at a general $1-$dimensional
family in $\n(r,d)$. Let $\C_0$ be a
general curve in $\n(r,d)$ that has
non-empty intersection with $\K$. First
we blow down the contracted components
of the fibres to get the family $\C$ where the marked points
may cross. The
family $\pi : S \to \C$ has two sections
$s_A,s_B : \C \to S$ corresponding two
preimages of the node, and there is a
map $\mu : S \to \proj$.
The two sections
cross at a point $P$. Let $Q = \pi(P)$ and
$F$ be the fibre over $Q$. We need to
show that $\mu_{|F}$ is not unramified at $P$.
This is easy to see. Let $\mu(s_A) = \mu(s_B) = G$
be the nodes of the image. For a general point $
t\in G$, the inverse image is a scheme of length
$2$ (in fact, it is a sum of two reduced points). Thus
$\mu^{-1}(\mu(P))$ is of length at least $2$. That means
the map $\mu_{|F}$ is not unramified at $P$. 
\epf

 We now give recursions to compute
characteristic numbers of rational cuspidal
curves, with condition on the cusp. Let
$p,q$ be two hyperplanes in $\proj.$ Consider
the moduli space $\X= \mbar_{0,S}(r,d)$ where
$S = \{A,B,P,Q\}$. Let $\n^S(r,d)$ be 
the closure in $\X$ of locus of maps $\gamma$
with $\gamma(A) = \gamma(B).$ We will
exploit the fact that if $\K$ is a boundary
divisor of $\X$ such that a general
map in $\K$ contracts a component containing
both $A$ and $B$, then $\n^S(r,d) \cap \K$
consists of rational cuspidal curves, with
the cusp at the image of $A$ (and $B$).
\begin{theorem}
Let $\d$ be a constraint. Let $k= \d(r+1).$
Let $\d_l$ be the same as $\d$ except that  $\d(0) -\d_l(0)=\d_l(r+1)-\d(r+1)= l$.
Let $m=\min(\d(0),r-\d(r+1))$. Let $\d'$ be the
same as $\d$ except that $\d'(2)= \d(2)+1$. Let
$\d''$ be the same as $\d$ except that
$\d''(r+1)= \d(r+1)+1$.
The following equality
holds if both sides are finite:
 \begin{eqnarray*}
d^2\#(\s(r,d),\d) &=& -\sum_{d_1+d_2=d}^{d_1,d_2>0}d_2^2\#(\n\R(r,d_1,d_2),\d)
 - \sum_{l=1}^m(^{\d(0)}_{ \ \ l})d^2\#(\s(r,d),\d_l) \\ 
&-& \#(\n(r,d),\d')  +\sum_{d_1+d_2=d}^{d_1,d_2>0}d_1d_2\#(\R\R_2(r,d_1,d_2),\d,k,0) \\
&+& 2d\#(\n(r,d),\d'').
 \end{eqnarray*}  
Furthermore, $\d_l$ is of lower rank than that of $\d.$
\end{theorem}
{\bf Warning :} if $\d(0) \neq 0$ then those summands above involving
reducible curves contain (twice) the case where the node
is mapped to a tangency hyperplane. Also, in computing
those summands, one needs to consider all possible splitting
of constraints over two components.

\bpf
Let $\wt{\d}$ be the same as $\d$ except
that $\wt{\d}$ contains the two hyperplanes
$p,q$. Let $\Y$ be the one-dimensional
family in $\n^S(r,d)$ of maps $\gamma$ that satisfy
$\d$ and also $\gamma(P) \in p$
and $\gamma(Q) \in q.$ We then intersect
$\Y$ with the two rationally equivalent
divisors:
$$(\{A,B\} \sep \{P,Q\}) \equiv (\{A,P\} \sep \{B, Q\})$$
First we intersect $\Y$ with the left
hand side of the equation. A general map $\gamma$ in
$(\{A,B\} \sep \{P,Q\}) $ has two-component
source curves. Let the component containing
$A,B$ be $C_1$, the other $C_2$. We consider 
cases:
\begin{itemize}
\item $\gamma_{|C_1}$ has degree $0$. By 
lemma $3.1$, we get $d^2\s(r,d) \cdot \d$.
The product is understood as intersecting
cycle classes representing constraint $\d$
on $\s(r,d)$. The factor $d^2$ comes from
two hyperplane conditions on $P$ and
$Q$. Since we are intersecting cycle 
classes, tangency to a hyperplane
means either ordinary tangency condition
where the cuspidal curve is tangent
to the hyperplane at a smooth point,
or it could mean the cusp lies on
the hyperplane. If we replane
$l$ ordinary tangency conditions
by $l$ hyperplane conditions through
the cusp, we get $\binom{\d(0)}{l}\#(\s(r,d),\d_l)$.
Thus the total contribution of this case
is 
\begin{eqnarray*}
d^2\#(\s(r,d),\d) + \sum_{l=1}^{\d(0)}d^2\binom{\d(0)}{l}d^2\#(\s(r,d),\d_l)\\ = d^2\#(\s(r,d),\d) + \sum_{l=1}^{m}d^2\binom{\d(0)}{l}d^2\#(\s(r,d),\d_l)
\end{eqnarray*}
\item $\gamma_{|C_2}$ has degree $0$. Two
hyperplane conditions on $P$ and $Q$ become
the top incident condition (that
explain why we need $\d'(2) = \d(2) +1$).
In this case we got $\#(\n(r,d),\d')$
\item $\gamma$ has positive degree on
each component. Let $d_i$ be the degree
of $\gamma$ on $C_i$. The total contribution
is 
$$\sum_{d_1+d_2=d}^{d_1,d_2>0}d_2^2\#(\n\R(r,d_1,d_2),\d)$$
The coefficient $d_2^2$ comes from
two hyperplane conditions on $P$ and 
$Q$ of the second component.
\end{itemize}
Now we intersection $\Y$ with the right
hand side. We consider $2$ cases:
\begin{itemize}
 \item $\gamma_{|C_1}$ has degree $0$ or 
$\gamma_{|C_2}$ has degree $0$. In each case,
the contribution is $d\#(\n(r,d),\d'')$. The total
contribution from both cases is $2d\#(\n(r,d),\d'').$
 \item $\gamma$ has positive degree on
each component. The contribution is
$$\sum_{d_1+d_2=d}^{d_1,d_2>0}d_1d_2\#(\R\R_2(r,d_1,d_2),\d,k,0) $$
\end{itemize}
Rearranging the terms we have the desired equation.
 \epf 

We have some flexibility in choosing the marked points
$A,B,P,Q$ to invoke the WDVV relation. For example,
we can impose (non-hyperplane) incident conditions
on $P$ and $Q$ and then invoke WDVV. In fact we have
a simpler and faster recursion this way with the
price of not being to compute all characteristic
numbers (for example, when all conditions are
tangency ). 

 Let $\d$ be a constraint such that $\d(0)=0$. Let $p,q$ be
two linear spaces in $\d$. Consider
the moduli space $\X= \mbar_{0,S}(r,d)$ where
$S = \{A,B,P,Q\}$. Let $\n^S(r,d)$ be 
the closure in $\X$ of locus of maps $\gamma$
with $\gamma(A) = \gamma(B)$. The following proposition produces
a different formula for incidence-only rational
cuspidal curve numbers.

\begin{proposition}
Let $k= \d(r+1).$
Let $\d'$ be derived from $\d$
by removing $p,q$ and add $p \cap q$.
Let $\d_p$ be derived from $\d$
by removing $p$ and add
$\codim(p)$ to $\d(r+1)$.
Let $\d_q$ be derived from $d$
by removing $q$ and add
$\codim(q)$ to $\d(r+1)$. Let
$\wt{\d}$ be derived from $\d$
by removing $p,q$. If $\gam$
is a constraint and $X$
is a set of linear spaces then $\gam^{(X)}$
is the constraint derived from
$\gam$ by adding linear spaces in $X$. 
The following equality
holds if both sides are finite:
\begin{eqnarray*}
 \#(\s(r,d),\d) &=& -\#(\n(r,d),\d') - \sum_{d_1+d_2=d}^{\gam_1\gam_2 = \wt{\d}}\binom{\wt{\d}} {\gam_1}\#(\n\R(r,d_1,d_2),\gam_1,\gam_2^{(p,q)},0) \\
&+& \sum_{d_1+d_2 = d}^{\gam_1\gam_2 = \wt{\d}}\binom{\wt{\d}}{\gam_1}\#(\R\R_2(r,d_1,d_2),\gam_1^{(p)},\gam_2^{(q)},k,0) \\
&+& \#(\n(r,d),\d_p) + \#(\n(r,d),\d_q).
\end{eqnarray*}
\end{proposition}
\bpf
The proof is identical to that of Theorem $3.2$.
\epf

 We implemented both recursions and confirm that
they gave the same numbers in case of incidence-only
constraints, which is good check of the method. 

 {\bf Formula for plane curves.}
  For plane curves with incident conditions, the formula simplify enough to get
 something computable by hands. We use the following
 standard notations. Let $C_d$ be the number of rational cuspidal plane curves,
 passing through $3d-2$ general points.
 Let $R_d$ be the number of rational plane curves passing
 through $3d-1$ points. Let $N_d,NL_d,NP_d$ be the
 number of rational plane curves with a choice of
 a node, where the node moves freely, on a line,
 on a point, respectively that pass through
 $3d-1,3d-2,3d-3$ points respectively.
 Now we restrict our formulas in Theorem $3.2$ and Proposition $3.3$
 to plane curves. Applying Proposition $3.3$ and Theorem $3.2$
 gives the following recursions for rational cuspidal plane
 curves with incident conditions.
 \begin{corollary}
 \begin{eqnarray*}
  C_d &=& 4NP_d + \sum_{i+j=d}\binom{3d-4}{3i-2}ij(ij-1)R_iR_{j} \\ 
 &-& \sum_{i=1}^{d-1}2\binom{3d-4}{3i-1}ijN_iR_{j}.
 \end{eqnarray*}
 \end{corollary}
 \begin{corollary}
 \begin{eqnarray*}
  d^2C_d &=& 4dNL_d + \sum_{i=1}^{d-1}\binom{3d-2}{3i-1}i^2j^2(ij-1)R_iR_j\\
  &-& 2N_d - \sum_{i+j=d}2\binom{3d-2}{3i-1}j^3iN_iR_j
 \end{eqnarray*}
 \end{corollary}
The formula for full chacracteristic numbers even
for plane curves is unfortunately still rather
involved (it may contain, for example, intersection
numbers on $Bl_{\D}(\mathbb{P}^2 \times \mathbb{P}^2),$ 
the blow up of product of projective planes along the
diagonal, see \cite{dn})  and  is best left
to a computer program.
\section{Numerical Examples}
 Let $C$ be the family of rational cuspidal
 curves. Let $C_f,C_b,C_s,C_l,C_p$ be the
 same families except that the cusp
 is required to lie on a $4-$space,
 $3-$space, a plane, a line, and a point respectively. The
 tables below list the characteristic numbers for
 such family. In particular, we recover all 
 previously known results, such as in \cite{luf},
 \cite{p}, \cite{hm}. In some of the tables,
we impose some point incident conditions
on the families to make the numbers small
enough to fit into the tables (tables
$10$ and $12$). All other conditions
are tangency and top incident conditions (incident with a codimension $2$ linear 
subspace). We were also able to recover
the results in table $5$ of \cite{zg} using
either of the recursions (in Theorem $3.2$
or Proposition $3.3$).

\begin{center}
 \begin{tabular} {| c | l | l | l |}
 \hline
$\#$ tang & $C$ & $C_l$ & $C_p$ \\
\hline
$ 0 $ &  $ 24 $ &  $ 12 $ &  $ 2 $   \\
\hline
$ 1 $ &  $ 60 $ &  $ 42 $ &  $ 8 $   \\
\hline
$ 2 $ &  $ 114 $ &  $ 96 $ &  $ 20 $   \\
\hline
$ 3 $ &  $ 168 $ &  $ 168 $ &  $ 38 $   \\
\hline
$ 4 $ &  $ 168 $ &  $ 186 $ &  $ 44 $  \\
\hline
$ 5 $ &  $ 114 $ &  $ 132 $ &  $ 32 $   \\
\hline
$ 6 $ &  $ 60 $ &  $ 72 $ &      \\
\hline
$ 7 $ &  $ 24 $ &    &      \\
\hline

 \end{tabular}

\end{center}
\begin{center}
  Table 1. Plane Cubics
\end{center}

\begin{center}
 \begin{tabular} {| c | l | l | l |}
 \hline
$\#$ tang & $C$ & $C_l$ & $C_p$ \\
\hline
$ 0 $ &  $ 2304 $ &  $ 864 $ &  $ 102 $   \\
\hline
$ 1 $ &  $ 6912 $ &  $ 2926 $ &  $ 364 $   \\
\hline
$ 2 $ &  $ 18486 $ &  $ 8552 $ &  $ 1112 $   \\
\hline
$ 3 $ &  $ 43668 $ &  $ 21816 $ &  $ 2964 $   \\
\hline
$ 4 $ &  $ 88560 $ &  $ 47284 $ &  $ 6700 $   \\
\hline
$ 5 $ &  $ 149364 $ &  $ 84284 $ &  $ 12392 $  \\
\hline
$ 6 $ &  $ 201132 $ &  $ 118296 $ &  $ 17904 $   \\
\hline
$ 7 $ &  $ 212976 $ &  $ 128992 $ &  $ 19912 $   \\
\hline
$ 8 $ &  $ 180288 $ &  $ 111776 $ &  $ 17444 $   \\
\hline
$ 9 $ &  $ 126720 $ &  $ 81324 $ &      \\
\hline
$ 10 $ &  $ 75924 $ &    &      \\
\hline

 \end{tabular}

\end{center}
\begin{center}
  Table 2. Plane quartics.
\end{center}

\begin{center} 
 \begin{tabular} {| c | l | l | l |} 
 \hline
$\#$ tang & $C$ & $C_l$ & $C_p$ \\
 \hline
$ 0 $ &  $ 435168 $ &  $ 130896 $ &  $ 12024 $  \\
\hline
$ 1 $ &  $ 1467792 $ &  $ 469112 $ &  $ 44272 $  \\
\hline
$ 2 $ &  $ 4592952 $ &  $ 1544416 $ &  $ 149504 $  \\
\hline
$ 3 $ &  $ 13240080 $ &  $ 4659408 $ &  $ 462744 $  \\
\hline
$ 4 $ &  $ 34794432 $ &  $ 12761768 $ &  $ 1300816 $  \\
\hline
$ 5 $ &  $ 82282248 $ &  $ 31325416 $ &  $ 3276944 $  \\
\hline
$ 6 $ &  $ 172272672 $ &  $ 67751352 $ &  $ 7264872 $  \\
\hline
$ 7 $ &  $ 313485192 $ &  $ 126634616 $ &  $ 13882384 $  \\
\hline
$ 8 $ &  $ 486730080 $ &  $ 200794048 $ &  $ 22428704 $  \\
\hline
$ 9 $ &  $ 637644672 $ &  $ 267383808 $ &  $ 30314016 $  \\
\hline
$ 10 $ &  $ 704860128 $ &  $ 299655536 $ &  $ 34336864 $  \\
\hline
$ 11 $ &  $ 664607952 $ &  $ 286715392 $ &  $ 33071072 $  \\
\hline
$ 12 $ &  $ 543805632 $ &  $ 239546016 $ &     \\
\hline
$ 13 $ &  $ 392798880 $ &    &     \\
\hline
 \end{tabular} 
 \end{center} 
 \begin{center}  
 Table $3$. Plane Quintics.
 \end{center}

\begin{center} 
 \begin{tabular} {| c | l | l | l |}
 \hline
$\#$ tang & $C$ & $C_l$ & $C_p$ \\ 
 \hline
$ 0 $ &  $ 156153600 $ &  $ 39223584 $ &  $ 2953656 $  \\
\hline
$ 1 $ &  $ 568978848 $ &  $ 148197528 $ &  $ 11344464 $  \\
\hline
$ 2 $ &  $ 1958182776 $ &  $ 526861728 $ &  $ 40992192 $  \\
\hline
$ 3 $ &  $ 6336731376 $ &  $ 1756900800 $ &  $ 138987648 $  \\
\hline
$ 4 $ &  $ 19174449024 $ &  $ 5467574880 $ &  $ 440024976 $  \\
\hline
$ 5 $ &  $ 53896263840 $ &  $ 15776422704 $ &  $ 1292243616 $  \\
\hline
$ 6 $ &  $ 139623413328 $ &  $ 41868541728 $ &  $ 3491056224 $  \\
\hline
$ 7 $ &  $ 330309167616 $ &  $ 101219247648 $ &  $ 8588991312 $  \\
\hline
$ 8 $ &  $ 706142918016 $ &  $ 220484782080 $ &  $ 19022973408 $  \\
\hline
$ 9 $ &  $ 1349047007424 $ &  $ 427799034720 $ &  $ 37472273856 $  \\
\hline
$ 10 $ &  $ 2279106096480 $ &  $ 731693677824 $ &  $ 64947355776 $  \\
\hline
$ 11 $ &  $ 3379250467008 $ &  $ 1095529768128 $ &  $ 98337086784 $  \\
\hline
$ 12 $ &  $ 4384937408256 $ &  $ 1433281143168 $ &  $ 129812346432 $  \\
\hline
$ 13 $ &  $ 4991823277056 $ &  $ 1644831322560 $ &  $ 149967485568 $  \\
\hline
$ 14 $ &  $ 5020701378624 $ &  $ 1670495736960 $ &  $ 152985851136 $  \\
\hline
$ 15 $ &  $ 4505776427520 $ &  $ 1519601842944 $ &     \\
\hline
$ 16 $ &  $ 3646864965888 $ &    &     \\
\hline
 \end{tabular} 
 \end{center} 
 \begin{center}  
 Table $4$. Plane sextics. 
 \end{center}

\begin{center} 
 \begin{tabular} {| c | l |  l | l | l |} 
 \hline 
 $\#$ tang & $C$ & $C_s$ & $C_l$ & $C_p$ \\
\hline
$ 0 $ &  $ 17760 $ &  $ 6592 $ &  $ 1168 $ &  $ 96 $  \\
\hline
$ 1 $ &  $ 31968 $ &  $ 14800 $ &  $ 2896 $ &  $ 264 $  \\
\hline
$ 2 $ &  $ 44304 $ &  $ 22336 $ &  $ 4592 $ &  $ 448 $  \\
\hline
$ 3 $ &  $ 49008 $ &  $ 25560 $ &  $ 5408 $ &  $ 556 $  \\
\hline
$ 4 $ &  $ 43104 $ &  $ 22864 $ &  $ 4952 $ &  $ 540 $  \\
\hline
$ 5 $ &  $ 30960 $ &  $ 16672 $ &  $ 3708 $ &  $ 436 $  \\
\hline
$ 6 $ &  $ 18888 $ &  $ 10380 $ &  $ 2376 $ &  $ 304 $  \\
\hline
$ 7 $ &  $ 10284 $ &  $ 5836 $ &  $ 1392 $ &  $ 208 $  \\
\hline
$ 8 $ &  $ 5088 $ &  $ 3040 $ &  $ 768 $ &     \\
\hline
$ 9 $ &  $ 2304 $ &  $ 1504 $ &    &     \\
\hline
$ 10 $ &  $ 960 $ &    &    &     \\
\hline
 \end{tabular} 
 \end{center} 
 \begin{center}  
 Table $5$. Cubics in $\mathbb{P}^3$. 
 \end{center}

\begin{center} 
 \begin{tabular} {| c | l |  l | l | l |} 
 \hline
$\#$ tang & $C$&  $C_s$ & $C_l$ & $C_p$ \\ \hline
$ 0 $ &  $ 170573760 $ &  $ 38051328 $ &  $ 4349376 $ &  $ 227088 $  \\
\hline
$ 1 $ &  $ 299241600 $ &  $ 72849984 $ &  $ 8741424 $ &  $ 479760 $  \\
\hline
$ 2 $ &  $ 476804928 $ &  $ 123397584 $ &  $ 15370848 $ &  $ 880992 $  \\
\hline
$ 3 $ &  $ 690153456 $ &  $ 186762096 $ &  $ 23996064 $ &  $ 1432240 $  \\
\hline
$ 4 $ &  $ 903026880 $ &  $ 252400320 $ &  $ 33289664 $ &  $ 2065784 $  \\
\hline
$ 5 $ &  $ 1062095040 $ &  $ 303769088 $ &  $ 40976856 $ &  $ 2643008 $  \\
\hline
$ 6 $ &  $ 1119351360 $ &  $ 325659864 $ &  $ 44849504 $ &  $ 3014080 $  \\
\hline
$ 7 $ &  $ 1059100728 $ &  $ 312649344 $ &  $ 43961088 $ &  $ 3092832 $  \\
\hline
$ 8 $ &  $ 906022656 $ &  $ 271494976 $ &  $ 39025664 $ &  $ 2892800 $  \\
\hline
$ 9 $ &  $ 708615360 $ &  $ 216184320 $ &  $ 31843584 $ &  $ 2507904 $  \\
\hline
$ 10 $ &  $ 513534720 $ &  $ 160405248 $ &  $ 24299136 $ &  $ 2053728 $  \\
\hline
$ 11 $ &  $ 349681920 $ &  $ 112854144 $ &  $ 17700384 $ &  $ 1654800 $  \\
\hline
$ 12 $ &  $ 226364928 $ &  $ 76461504 $ &  $ 12539328 $ &     \\
\hline
$ 13 $ &  $ 140561856 $ &  $ 50638656 $ &    &     \\
\hline
$ 14 $ &  $ 84248640 $ &    &    &     \\
\hline
 \end{tabular} 
 \end{center} 
 \begin{center}  
 Table $6$. Quartics in $\mathbb{P}^3$. 
 \end{center}

\begin{center} 
 \begin{tabular} {| c | l |  l | l | l |} 
 \hline
$\#$ tang & $C$&  $C_s$ & $C_l$ & $C_p$ \\ \hline
$ 0 $ &  $ 3367116891648 $ &  $ 550184948736 $ &  $ 46665372672 $ &  $ 1785262368 $  \\
\hline
$ 1 $ &  $ 6238276331520 $ &  $ 1067735961600 $ &  $ 93212510688 $ &  $ 3676278528 $  \\
\hline
$ 2 $ &  $ 10907520952320 $ &  $ 1939988681376 $ &  $ 173756646144 $ &  $ 7055672256 $  \\
\hline
$ 3 $ &  $ 17959194667872 $ &  $ 3299107194624 $ &  $ 302481758400 $ &  $ 12639091008 $  \\
\hline
$ 4 $ &  $ 27757496169984 $ &  $ 5240093299584 $ &  $ 490925857152 $ &  $ 21103670976 $  \\
\hline
$ 5 $ &  $ 40132781267328 $ &  $ 7751902200576 $ &  $ 740922699072 $ &  $ 32768348448 $  \\
\hline
$ 6 $ &  $ 54102237425280 $ &  $ 10651963788480 $ &  $ 1037273078208 $ &  $ 47214082368 $  \\
\hline
$ 7 $ &  $ 67828147429440 $ &  $ 13570111226976 $ &  $ 1344945050304 $ &  $ 63061854976 $  \\
\hline
$ 8 $ &  $ 78978275227392 $ &  $ 16020342530304 $ &  $ 1615192092160 $ &  $ 78141404928 $  \\
\hline
$ 9 $ &  $ 85438440716544 $ &  $ 17550106827264 $ &  $ 1800060152832 $ &  $ 90077596416 $  \\
\hline
$ 10 $ &  $ 86056052126208 $ &  $ 17898091977216 $ &  $ 1868700617472 $ &  $ 97050126528 $  \\
\hline
$ 11 $ &  $ 81016897731840 $ &  $ 17076676508160 $ &  $ 1817043169728 $ &  $ 98356055808 $  \\
\hline
$ 12 $ &  $ 71671571228160 $ &  $ 15340937774976 $ &  $ 1666394974464 $ &  $ 94521559296 $  \\
\hline
$ 13 $ &  $ 59960396118912 $ &  $ 13074466454784 $ &  $ 1453197160320 $ &  $ 86990259648 $  \\
\hline
$ 14 $ &  $ 47770968040704 $ &  $ 10659479487744 $ &  $ 1216151757312 $ &  $ 77589435456 $  \\
\hline
$ 15 $ &  $ 36506631132288 $ &  $ 8387937598656 $ &  $ 986787943872 $ &  $ 68371627008 $  \\
\hline
$ 16 $ &  $ 26947730227200 $ &  $ 6429428144640 $ &  $ 784978434048 $ &     \\
\hline
$ 17 $ &  $ 19339799508480 $ &  $ 4848964862976 $ &    &     \\
\hline
$ 18 $ &  $ 13574889925632 $ &    &    &     \\
\hline
 \end{tabular} 
 \end{center} 
 \begin{center}  
 Table $7$. Quintics in $\mathbb{P}^3$. 
 \end{center}
\begin{center} 
 \begin{tabular} {| c | l |  l | l | l | l |} 
 \hline
$\#$ tang & $C$& $C_b$ & $C_s$ & $C_l$ & $C_p$ \\ \hline
$ 0 $ &  $ 8004000 $ &  $ 2444280 $ &  $ 435550 $ &  $ 48660 $ &  $ 2746 $  \\
\hline
$ 1 $ &  $ 11533080 $ &  $ 4107630 $ &  $ 790900 $ &  $ 95094 $ &  $ 5824 $  \\
\hline
$ 2 $ &  $ 13193310 $ &  $ 4934340 $ &  $ 981226 $ &  $ 122752 $ &  $ 7916 $  \\
\hline
$ 3 $ &  $ 12522180 $ &  $ 4727094 $ &  $ 954432 $ &  $ 122692 $ &  $ 8254 $  \\
\hline
$ 4 $ &  $ 10083450 $ &  $ 3799648 $ &  $ 776660 $ &  $ 102828 $ &  $ 7290 $  \\
\hline
$ 5 $ &  $ 7066080 $ &  $ 2660080 $ &  $ 552170 $ &  $ 75800 $ &  $ 5726 $  \\
\hline
$ 6 $ &  $ 4422480 $ &  $ 1673790 $ &  $ 354420 $ &  $ 50714 $ &  $ 4124 $  \\
\hline
$ 7 $ &  $ 2537070 $ &  $ 973580 $ &  $ 211186 $ &  $ 31692 $ &  $ 2816 $  \\
\hline
$ 8 $ &  $ 1362780 $ &  $ 534914 $ &  $ 119572 $ &  $ 19032 $ &  $ 1904 $  \\
\hline
$ 9 $ &  $ 692850 $ &  $ 280908 $ &  $ 65080 $ &  $ 11088 $ &  $ 1280 $  \\
\hline
$ 10 $ &  $ 335460 $ &  $ 142080 $ &  $ 34240 $ &  $ 6240 $ &     \\
\hline
$ 11 $ &  $ 155280 $ &  $ 69600 $ &  $ 17440 $ &    &     \\
\hline
$ 12 $ &  $ 68880 $ &  $ 33120 $ &    &    &     \\
\hline
$ 13 $ &  $ 29280 $ &    &    &    &     \\
\hline
 \end{tabular} 
 \end{center} 
 \begin{center}  
 Table $8$. Cubics in $\mathbb{P}^4$. 
 \end{center}

\begin{center} 
 \begin{tabular} {| c | l |  l | l | l | l |} 
 \hline
$\#$ tang & $C$& $C_b$ & $C_s$ & $C_l$ & $C_p$ \\ \hline
$ 0 $ &  $ 2689538273280 $ &  $ 485362519880 $ &  $ 54851631000 $ &  $ 3880179030 $ &  $ 135655104 $  \\
\hline
$ 1 $ &  $ 3724396027560 $ &  $ 710716985200 $ &  $ 83317313210 $ &  $ 6115563684 $ &  $ 222942852 $  \\
\hline
$ 2 $ &  $ 4732678206600 $ &  $ 937159830750 $ &  $ 112867583584 $ &  $ 8537849664 $ &  $ 322772136 $  \\
\hline
$ 3 $ &  $ 5545209427650 $ &  $ 1126146131924 $ &  $ 138560376836 $ &  $ 10764060168 $ &  $ 421095648 $  \\
\hline
$ 4 $ &  $ 6006642647520 $ &  $ 1241407662880 $ &  $ 155502195280 $ &  $ 12384357920 $ &  $ 501209960 $  \\
\hline
$ 5 $ &  $ 6029598965220 $ &  $ 1261851548400 $ &  $ 160600815680 $ &  $ 13107423520 $ &  $ 549512524 $  \\
\hline
$ 6 $ &  $ 5627161190520 $ &  $ 1189150782880 $ &  $ 153667840720 $ &  $ 12862780284 $ &  $ 560131552 $  \\
\hline
$ 7 $ &  $ 4904324736960 $ &  $ 1045491288800 $ &  $ 137211306524 $ &  $ 11798753624 $ &  $ 535636616 $  \\
\hline
$ 8 $ &  $ 4014398317200 $ &  $ 863628780684 $ &  $ 115219223936 $ &  $ 10199296568 $ &  $ 484737728 $  \\
\hline
$ 9 $ &  $ 3106478617740 $ &  $ 675426920360 $ &  $ 91724366160 $ &  $ 8377299120 $ &  $ 418707360 $  \\
\hline
$ 10 $ &  $ 2288765734560 $ &  $ 504110473120 $ &  $ 69796390160 $ &  $ 6591766560 $ &  $ 348145344 $  \\
\hline
$ 11 $ &  $ 1617094353480 $ &  $ 361922529360 $ &  $ 51181998880 $ &  $ 5009865984 $ &  $ 281127552 $  \\
\hline
$ 12 $ &  $ 1103132481120 $ &  $ 251838754400 $ &  $ 36454172864 $ &  $ 3708038784 $ &  $ 222727296 $  \\
\hline
$ 13 $ &  $ 730910337600 $ &  $ 170972796544 $ &  $ 25393735936 $ &  $ 2692309568 $ &  $ 174850208 $  \\
\hline
$ 14 $ &  $ 472615004160 $ &  $ 113850253440 $ &  $ 17390191040 $ &  $ 1925281120 $ &  $ 136210960 $  \\
\hline
$ 15 $ &  $ 299278396800 $ &  $ 74653047040 $ &  $ 11745665440 $ &  $ 1355930400 $ &     \\
\hline
$ 16 $ &  $ 186019069440 $ &  $ 48323296640 $ &  $ 7834203440 $ &    &     \\
\hline
$ 17 $ &  $ 113612974080 $ &  $ 30909046800 $ &    &    &     \\
\hline
$ 18 $ &  $ 68177775600 $ &    &    &    &     \\
\hline
 \end{tabular} 
 \end{center} 
 \begin{center}  
 Table $9$. Quartics in $\mathbb{P}^4.$ 
 \end{center}
\begin{center} 
 \begin{tabular} {| c | l |  l | l | l | l |} 
 \hline
$\#$ tang & $C$& $C_b$ & $C_s$ & $C_l$ & $C_p$ \\ \hline
$ 0 $ &  $ 326137674360 $ &  $ 49962812760 $ &  $ 5012100872 $ &  $ 341811654 $ &  $ 12360978 $  \\
\hline
$ 1 $ &  $ 525454848600 $ &  $ 84239622360 $ &  $ 8769540462 $ &  $ 628334598 $ &  $ 24199100 $  \\
\hline
$ 2 $ &  $ 810363195432 $ &  $ 135076427046 $ &  $ 14543938282 $ &  $ 1092843240 $ &  $ 44814208 $  \\
\hline
$ 3 $ &  $ 1197614346510 $ &  $ 206664134118 $ &  $ 22970962524 $ &  $ 1809906216 $ &  $ 79191906 $  \\
\hline
$ 4 $ &  $ 1697274024738 $ &  $ 302362862376 $ &  $ 34656015928 $ &  $ 2865776848 $ &  $ 134313322 $  \\
\hline
$ 5 $ &  $ 2308283286660 $ &  $ 423823209200 $ &  $ 50062418802 $ &  $ 4351800704 $ &  $ 219711358 $  \\
\hline
$ 6 $ &  $ 3015524277240 $ &  $ 570330420630 $ &  $ 69407614892 $ &  $ 6356646798 $ &  $ 348417840 $  \\
\hline
$ 7 $ &  $ 3789887203782 $ &  $ 738694207496 $ &  $ 92608706262 $ &  $ 8960124980 $ &  $ 538817608 $  \\
\hline
$ 8 $ &  $ 4591908391500 $ &  $ 923850975318 $ &  $ 119302033904 $ &  $ 12228114496 $ &  $ 817997616 $  \\
\hline
$ 9 $ &  $ 5378099143998 $ &  $ 1120077375516 $ &  $ 148922757248 $ &  $ 16197194208 $ &  $ 1226358432 $  \\
\hline
$ 10 $ &  $ 6107875465560 $ &  $ 1322485952040 $ &  $ 180803857632 $ &  $ 20819978304 $ &  $ 1819709888 $  \\
\hline
$ 11 $ &  $ 6748291059480 $ &  $ 1528252628400 $ &  $ 214285519072 $ &  $ 25840126848 $ &     \\
\hline
$ 12 $ &  $ 7273256585232 $ &  $ 1736403696096 $ &  $ 248858380992 $ &    &     \\
\hline
$ 13 $ &  $ 7655587765440 $ &  $ 1944122007168 $ &    &    &     \\
\hline
$ 14 $ &  $ 7857466535808 $ &    &    &    &     \\
\hline
 \end{tabular} 
 \end{center} 
 \begin{center}  
 Table $10.$ Quintics in $\mathbb{P}^5$ passing through $3$ points. 
 \end{center}

\begin{center} 
 \begin{tabular} {| c | l |  l | l | l | l | l |} 
 \hline
$\#$ tang & $C$& $C_f$ & $C_b$ & $C_s$ & $C_l$ & $C_p$ \\ \hline
$ 0 $ &  $ 2872888704 $ &  $ 764380848 $ &  $ 134764896 $ &  $ 16936134 $ &  $ 1424016 $ &  $ 61930 $  \\
\hline
$ 1 $ &  $ 3520649664 $ &  $ 1046482248 $ &  $ 196268970 $ &  $ 26096280 $ &  $ 2330142 $ &  $ 108400 $  \\
\hline
$ 2 $ &  $ 3534202224 $ &  $ 1084974798 $ &  $ 208339392 $ &  $ 28472226 $ &  $ 2630496 $ &  $ 127672 $  \\
\hline
$ 3 $ &  $ 3038206194 $ &  $ 935625864 $ &  $ 181474518 $ &  $ 25249632 $ &  $ 2396376 $ &  $ 120756 $  \\
\hline
$ 4 $ &  $ 2297295216 $ &  $ 703602714 $ &  $ 137487936 $ &  $ 19484376 $ &  $ 1906920 $ &  $ 100592 $  \\
\hline
$ 5 $ &  $ 1562512590 $ &  $ 476145840 $ &  $ 93931320 $ &  $ 13610844 $ &  $ 1381596 $ &  $ 76886 $  \\
\hline
$ 6 $ &  $ 975034080 $ &  $ 297003360 $ &  $ 59370960 $ &  $ 8831220 $ &  $ 934260 $ &  $ 55260 $  \\
\hline
$ 7 $ &  $ 568215648 $ &  $ 174117384 $ &  $ 35401716 $ &  $ 5424480 $ &  $ 601068 $ &  $ 38116 $  \\
\hline
$ 8 $ &  $ 313900272 $ &  $ 97387668 $ &  $ 20211060 $ &  $ 3203004 $ &  $ 374160 $ &  $ 25744 $  \\
\hline
$ 9 $ &  $ 166136292 $ &  $ 52496400 $ &  $ 11160396 $ &  $ 1837716 $ &  $ 228096 $ &  $ 17296 $  \\
\hline
$ 10 $ &  $ 84772620 $ &  $ 27438780 $ &  $ 5994360 $ &  $ 1029408 $ &  $ 136608 $ &  $ 11584 $  \\
\hline
$ 11 $ &  $ 41869260 $ &  $ 13961244 $ &  $ 3141360 $ &  $ 563712 $ &  $ 80256 $ &  $ 7744 $  \\
\hline
$ 12 $ &  $ 20069856 $ &  $ 6933936 $ &  $ 1608768 $ &  $ 301632 $ &  $ 46080 $ &     \\
\hline
$ 13 $ &  $ 9353952 $ &  $ 3367920 $ &  $ 805824 $ &  $ 157632 $ &    &     \\
\hline
$ 14 $ &  $ 4243680 $ &  $ 1602048 $ &  $ 395136 $ &    &    &     \\
\hline
$ 15 $ &  $ 1874880 $ &  $ 747264 $ &    &    &    &     \\
\hline
$ 16 $ &  $ 806400 $ &    &    &    &    &     \\
\hline
 \end{tabular} 
 \end{center} 
 \begin{center}  
 Table $11.$ Cubics in $\mathbb{P}^5$. 
 \end{center}

\begin{center} 
 \begin{tabular} {| c | l |  l | l | l | l | l | } 
 \hline
$\#$ tang & $C$& $C_f$ & $C_b$ & $C_s$ & $C_l$ & $C_p$ \\ \hline
$ 0 $ &  $ 11321640 $ &  $ 2552292 $ &  $ 401422 $ &  $ 51972 $ &  $ 5178 $ &  $ 306 $  \\
\hline
$ 1 $ &  $ 17941716 $ &  $ 4364030 $ &  $ 726024 $ &  $ 101982 $ &  $ 11310 $ &  $ 768 $  \\
\hline
$ 2 $ &  $ 26377710 $ &  $ 6787492 $ &  $ 1181426 $ &  $ 178170 $ &  $ 21828 $ &  $ 1704 $  \\
\hline
$ 3 $ &  $ 36390888 $ &  $ 9813894 $ &  $ 1774470 $ &  $ 285216 $ &  $ 38472 $ &  $ 3492 $  \\
\hline
$ 4 $ &  $ 47498634 $ &  $ 13374146 $ &  $ 2494892 $ &  $ 423696 $ &  $ 62632 $ &  $ 6742 $  \\
\hline
$ 5 $ &  $ 59064222 $ &  $ 17368048 $ &  $ 3317488 $ &  $ 588180 $ &  $ 94208 $ &  $ 12416 $  \\
\hline
$ 6 $ &  $ 70392240 $ &  $ 21701328 $ &  $ 4211904 $ &  $ 768198 $ &  $ 129420 $ &     \\
\hline
$ 7 $ &  $ 80723952 $ &  $ 26295088 $ &  $ 5154310 $ &  $ 957324 $ &    &     \\
\hline
$ 8 $ &  $ 89121888 $ &  $ 31052678 $ &  $ 6119300 $ &    &    &     \\
\hline
$ 9 $ &  $ 94411494 $ &  $ 35864316 $ &    &    &    &     \\
\hline
$ 10 $ &  $ 95150484 $ &    &    &    &    &     \\
\hline
 \end{tabular} 
 \end{center} 
 \begin{center}  
 Table $12$. Quartics in $\mathbb{P}^5$ passing through $3$ points. 
 \end{center}

\end{document}